%% file: main.tex
\documentclass[a4paper]{article}
\usepackage{amssymb}
\usepackage{amsmath}
\usepackage{epsfig}
\usepackage{xspace}
\input latex.def

\newcommand{\card}[1]{\lvert {#1} \rvert}

\newcommand{\polymake}{\texttt{polymake}\xspace}
\newcommand{\javaview}{\texttt{javaview}\xspace}

\newcommand{\ignore}[1]{}
\newcommand{\verts}[1]{\text{\rm vert}({#1})}
\newcommand{\fingr}[1]{\Gamma_{#1}}

\title{Vertex-Facet Incidences of Unbounded Polyhedra}
\author{Michael Joswig \and
Volker Kaibel$^*$ \and
Marc~E.~Pfetsch \and
G\"unter M. Ziegler\thanks{%
Supported by a DFG Gerhard-Hess-Forschungsf\"orderungspreis (Zi 475/2-3).}\\[2mm]
Dept.\ Mathematics, MA~7-1\\ 
Technische Universit\"at Berlin\\
10623 Berlin, Germany\\
\{\texttt{joswig,kaibel,pfetsch,ziegler}\}\texttt{@math.tu-berlin.de}}

\begin{document}
\maketitle

\begin{abstract}
How much of the combinatorial structure of a pointed polyhedron
is contained in its vertex-facet incidences?
Not too much, in general, as we demonstrate by examples.
However, one can tell from the incidence data whether the polyhedron
is bounded. In the case of a polyhedron that is simple and
``simplicial,'' i.e., a $d$-dimensional polyhedron that has $d$ facets
through each vertex and $d$~vertices on each facet, we derive from
the structure of the vertex-facet incidence matrix that the polyhedron
is necessarily bounded. In particular, this yields a characterization
of those polyhedra that have circulants as vertex-facet incidence
matrices. 
\end{abstract}

\input intro.tex
\input definitions.tex
\input reconstructing.tex
\input unbounded.tex
\input simpsimp.tex

\bibliographystyle{amsplain}
\bibliography{polyhedra}

\end{document}

%% file: intro.tex
\section{Introduction}
\label{sec:intro}

Every (proper) face of a polytope (i.e., a bounded convex
polyhedron) is the convex hull of the vertices it contains, and it is
also the intersection of the facets that contain it. Thus, the
combinatorial structure of a polytope (i.e., its face lattice) is
entirely determined by its (matrix of) vertex-facet incidences.  Such
a vertex-facet incidence matrix is a useful encoding of the
combinatorial structure of a polytope. The software package {\tt
  polymake} \cite{pm00,GJ00}, for instance, represents this matrix
rather compactly, in a section called {\tt VERTICES\_IN\_FACETS},
while the face lattice of a polytope is not stored, but generated ``on
demand'' only if this is really necessary, because typically the
entire face lattice is ``much too large.''

But how about not necessarily bounded convex polyhedra?  The
combinatorics of unbounded polyhedra has received only little
attention up to now (for some exceptions see Klee~\cite{Klee74},
Billera \& Lee~\cite{BL81}, Barnette, Kleinschmidt \& Lee \cite{BKL},
and Lee~\cite{Lee84}).  One can, of course, reduce the study of
geometrically given unbounded polyhedra to the situation of ``a
polytope with a distinguished face (at infinity).''  But what if only
the combinatorics of vertices versus facets is given, and not any data
about the situation ``at infinity?''  In other words, how much can
really be said/detected/reconstructed if only a matrix of the
vertex-facet incidences is given?

As one observes easily from the example of polyhedral cones, in
general the combinatorial structure of an unbounded polyhedron is not
determined by its vertex-facet incidences. A $d$-dimensional cone 
may have any possible combinatorial structure of a
$(d-1)$-dimensional polytope (via \emph{homogenization}); but from its
vertex-facet incidences one can read off only its number of
facets.  The point is that, for unbounded polyhedra the
combinatorial information is based  not only on the vertex-facet
incidences, but also on the incidences of  extremal rays  and facets.
For cones, nearly the entire information is contained in the latter
incidences. The lattice-theoretic reason for such ambiguities is that
the face lattice of an unbounded polyhedron is only
co-atomic, but not atomic.

One might, however, suspect that cones are (extreme) examples of
rather exotic unbounded polyhedra for which one obviously does not
have any chance to reconstruct the combinatorial structure from their
vertex-facet incidences, while this might be possible for all
``reasonable'' polyhedra. For instance, a cone is a quite degenerate
polyhedron with respect to several criteria: (i) all its facets have
the same set of vertices, (ii) its set of vertices does not have the
same dimension as the whole polyhedron, and (iii) it does not have any
bounded facet.  However, the first main point of this paper (in
Section~\ref{sec:reconstruction}) is the construction of more
convincing examples of unbounded polyhedra whose face-lattices cannot
be reconstructed from their vertex-facet incidence matrices; they have
the property that the sets of vertices of facets are distinct, and
they even form an anti-chain in the Boolean lattice (a
\emph{clutter}); they have bounded facets, and their sets of vertices
are full-dimensional.


The second main result (in Section~\ref{sec:unbounded}) will be that
one can, however, detect from the vertex-facet matrix whether the
polyhedron under consideration is bounded or not.

Thirdly (in Section~\ref{sec:simpsimp}), we discuss the ``unbounded
version'' of a very basic lemma about polytopes. Indeed, Exercise 0.1
of \cite{GMZ} asks one to prove that any $d$-polytope that is both
simplicial (every facet has $d$ vertices) and simple (every vertex is
on $d$ facets) must either be a simplex, or a polygon ($d=2$).  But
how about unbounded polyhedra? We prove that a polyhedron that is both
simple and simplicial (with the definitions as given here) cannot be
unbounded.  As a byproduct, we obtain a characterization of those
polyhedra that have circulant vertex-facet incidence matrices.

\smallskip

In particular, this paper answers a series of questions that arose in 
Amaldi, Pfetsch, and Trotter~\cite{AmaPT00}, where the structure of certain
independence systems is related to the combinatorics of 
(possibly unbounded) polyhedra.


%% file: definitions.tex
\section{Basic Facts}
\label{sec:definitions}

Let $P$ be a $d$-polyhedron (i.e., the intersection of a finite number 
of affine halfspaces with $\dim(P)=d$) with $m$ facets and
$n$ vertices. We will always assume that $P$ is pointed (i.e., it has
at least one vertex) and that $d\geq 1$. In particular, these
conditions imply $n\geq 1$ and $m\geq d\geq 1$.  For the basic
definitions and facts of polyhedral theory we refer to~\cite{GMZ}.

A $0/1$-matrix $A=(a_{fv})\in\{0,1\}^{m\times n}$ is a
\emph{vertex-facet incidence matrix} of $P$ if the vertices and facets
of $P$ can be numbered by $\{1,\dots,n\}$ and $\{1,\dots,m\}$,
respectively, such that $a_{fv}=1$ if and only if the vertex with
number $v$ is contained in the facet with number $f$.

\begin{figure}[htbp]
  \begin{center}
    \epsfig{file=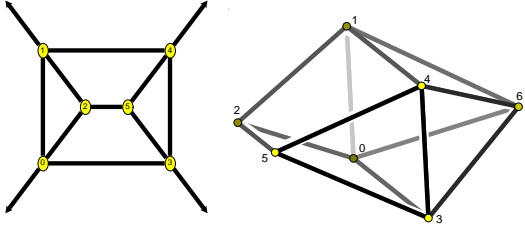,width=3cm}
    \epsfig{file=P_.eps,width=4cm}
    \caption{Left: $1$-skeleton (i.e., $0$- and $1$-dimensional faces) of a $3$-polyhedron~$P$.  The arrows indicate 
      extremal rays, which are assumed to be parallel.  Right: Projectively
      transformed into~$\overline{P}$.  The far face~$F_\infty$ is the
      vertex~$6$.}
    \label{fig:P}
  \end{center}
\end{figure}

By $\projclo{P}$ we denote any \emph{polytope} which is projectively
equivalent to $P$. If $P$ is unbounded, then there is a unique maximal
element $\farface$ (the \emph{far face}) among the faces of
$\projclo{P}$ that are not images of faces of $P$ under the projective
transformation mapping $P$ to $\projclo{P}$. If $P$ is bounded, then
we define $\farface=\varnothing$. Figure~\ref{fig:P} illustrates a
three-dimensional example.

We denote by $\facposet{P}$ the \emph{face poset} of $P$, i.e., the
set of non-trivial faces of $P$ (excluding $\varnothing$ and $P$
itself), ordered by inclusion. The face poset $\facposet{P}$ arises
from the face poset $\facposetprojclo{P}$ by removing the far face
$\farface$ (and all its faces).  While $\facposetprojclo{P}$ is
independent of the actual choice of $\projclo{P}$, in general it
depends on the geometry of $P$, not only on its combinatorial
structure.

\begin{figure}[htbp]
  \begin{center}
    \epsfig{file=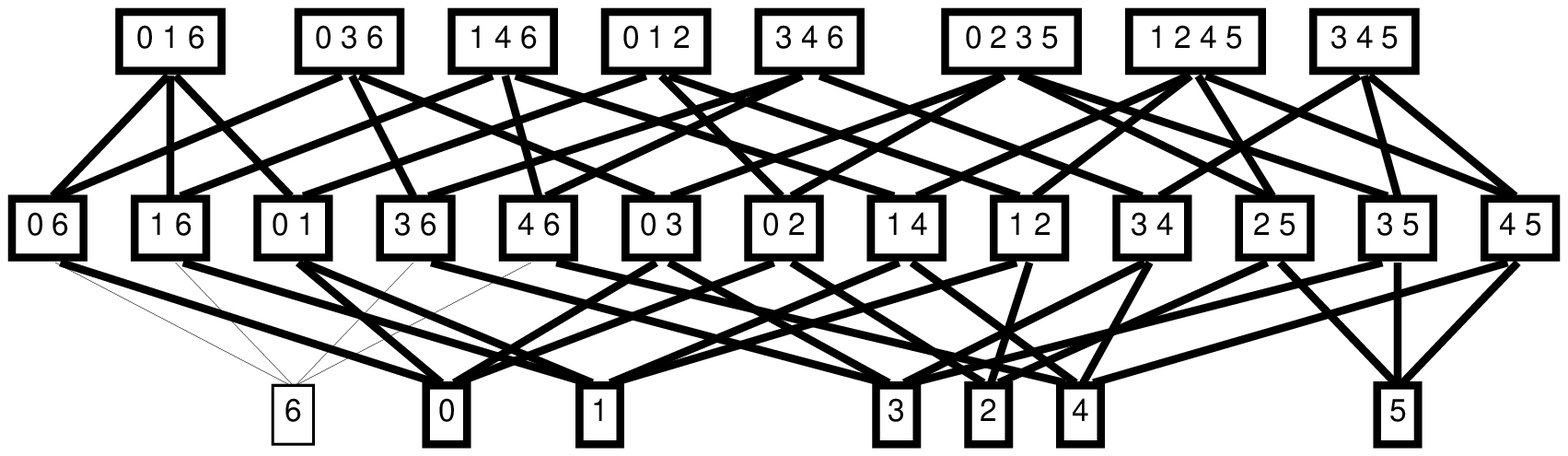,height=1.8cm}\quad
    \epsfig{file=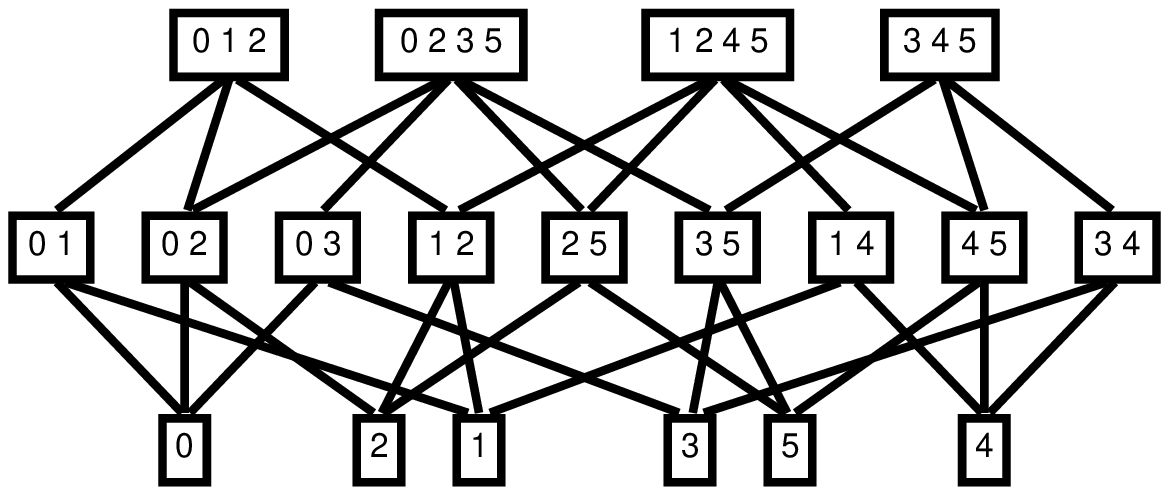,height=1.8cm}
    \caption{Left: Face poset of $\facposetprojclo{P}$ for the example
      given in Figure~\ref{fig:P}, where the solid part is
      $\facposet{P}$. Right: The poset $\facvertposet{P}$. 
      In general, $\facvertposet{P}$ is not graded (although it is in
      our example).}
    \label{fig:Plattices}
  \end{center}
\end{figure}

The poset $\facvertposet{P}=\SetOf{\verts{F}}{F\text{ non-trivial
    face of }P}$ (where $\verts{F}$ is the set of vertices of $F$)
will play an important role.  It can be computed from any vertex-facet
incidence matrix $A\in\{0,1\}^{m\times n}$ of $P$, since it is the set
of all non-empty intersections of subsets of $\{1,\dots,n\}$ defined
by subsets of the rows of $A$.  Figure~\ref{fig:Plattices} shows the
three posets $\facposet{P}$, $\facposetprojclo{P}$, and
$\facvertposet{P}$ for the example given in Figure~\ref{fig:P}.


Let the \emph{graph} $\fingr{P}$ of $P$ be the graph on the vertices
of $P$ defined by the bounded one-dimensional faces of $P$ (the
\emph{edges}), i.e., $\fingr{P}$ is the subgraph of the graph of
$\projclo{P}$ that is induced by those vertices of $\projclo{P}$ that
are not contained in $\farface$.  Two vertices of $P$ are connected by
an edge of $P$ if and only if there is a face of $P$ which contains
exactly these two vertices.  Moreover, we can compute
$\facvertposet{P}$ from any vertex-facet incidence matrix of $P$. In
particular, we can find $\fingr{P}$ from the vertex-facet incidences
of $P$.

We will use the following fact, which is a consequence of the
correctness of the Simplex-Algorithm for Linear Programming.

\begin{lem}\label{lem:gamma_connected}
  For every polyhedron $P$, the graph $\fingr{P}$ is connected.
  Moreover, all faces of~$P$ induce connected subgraphs of~$\fingr{P}$.
\end{lem}

Let~$P$ be a pointed $d$-polyhedron ($d\geq 1$).  Then, $P$ is called
\emph{simple} if every vertex of $P$ is contained in precisely $d$
facets (or, equivalently, if precisely $d$ edges and extremal rays are
incident to each vertex), and~$P$ is called \emph{simplicial} if every
facet of $P$ has precisely $d$ vertices. These notions generalize the
well-known notions \emph{simple} and \emph{simplicial} for
polytopes.  While this generalization is standard for simple
polyhedra, it is not common for simplicial polyhedra. Thus, it seems
to be worth to mention that simplicial unbounded polyhedra form a
non-trivial class of polyhedra. For instance, by a modification of the
construction of a prism, one easily sees that every simplicial
$d$-polytope can occur as the far face of a simplicial unbounded
$(d+1)$-polyhedron.


%% file: reconstructing.tex
\section{Reconstructing Polyhedra from Vertex-Facet Incidences}
\label{sec:reconstruction}

In this section, we consider conditions under which it is possible to
compute $\facposet{P}$ from the vertex-facet incidences of an
(unbounded) $d$-polyhedron $P$. Obviously, given any vertex-facet
incidence matrix of a pointed $d$-polyhedron $P$ it is easy to decide
whether $d\in\{1,2\}$.  Furthermore, if $d\in\{1,2\}$, one can immediately
read off $\facposet{P}$ from the vertex-facet incidences. Thus, for
the rest of this section we restrict our attention to $d$-polyhedra
with $d\geq 3$.

The example of cones shows that reconstructing $\facposet{P}$ from the
vertex-facet incidences of a $d$-polyhedron $P$ with $d\geq 4$ is
impossible in general, even if additionally the dimension $d$ is
specified. Furthermore, the same example demonstrates that it is, in
general, impossible to detect the dimension of a $d$-polyhedron from
its vertex-facet incidences for $d\geq 3$. However, for $d=3$ these
dimensional ambiguities occur for cones only.

\begin{prp}
\label{prp:detectDimThree}
  Given a vertex-facet incidence matrix of a $d$-polyhedron~$P$ with
  $d\geq 3$, it is
  possible to decide whether $d=3$ or $d\geq 4$, unless 
  $P$ is a cone with more than three facets.
\end{prp}

\begin{proof}
  If $P$ is a cone with three facets (i.e., $n=1$ and $m=3$) then
  clearly $d=3$ holds. If~$P$ is not a cone, then it
  must have at least two vertices.  Thus (by
  Lemma~\ref{lem:gamma_connected}) $P$ has at least one edge (which we
  can tell from the vertex-facet incidences of $P$). This edge is
  contained in precisely two facets of $P$ if $d=3$; otherwise, it
  is contained in more than two facets.
\end{proof}

In dimensions larger than three, cones are not the only polyhedra for
which one cannot tell the dimension from the vertex-facet incidences.
For instance, let $Q$ be some $d'$-polytope and let $C$ be a
$d''$-dimensional polyhedral cone with $m\geq 4$ facets. Then
$P=Q\times C$ will be a $(d'+d'')$-dimensional polyhedron whose
vertex-facet incidences only depend on $Q$ and $m$, while its
dimension can be any number between $d'+3$ and $d'+m$. In
particular, dimensional ambiguities already occur for $4$-polyhedra
not being cones.

However, the cartesian products constructed above are also
``cone-like'' in the sense that they do not have any bounded facet.

\begin{prp}
\label{prp:bndFacetDim}
  Given a vertex-facet incidence matrix of a $d$-polyhedron $P$ that
  has a bounded facet, one can determine $d$. Furthermore, one can
  decide from the vertex-facet incidences of $P$ whether it has a
  bounded facet or not.
\end{prp}

\begin{proof}
  If $P$ has a bounded facet, then the maximum length of a chain in
  $\facvertposet{P}$ is $d-1$, thus one can compute $d$ from
  $\facvertposet{P}$ in this case. Corollary~\ref{cor:facetDecide}
  proves the second statement of the proposition.
\end{proof}

\begin{figure}[htbp]
  \begin{center}
    \epsfig{figure=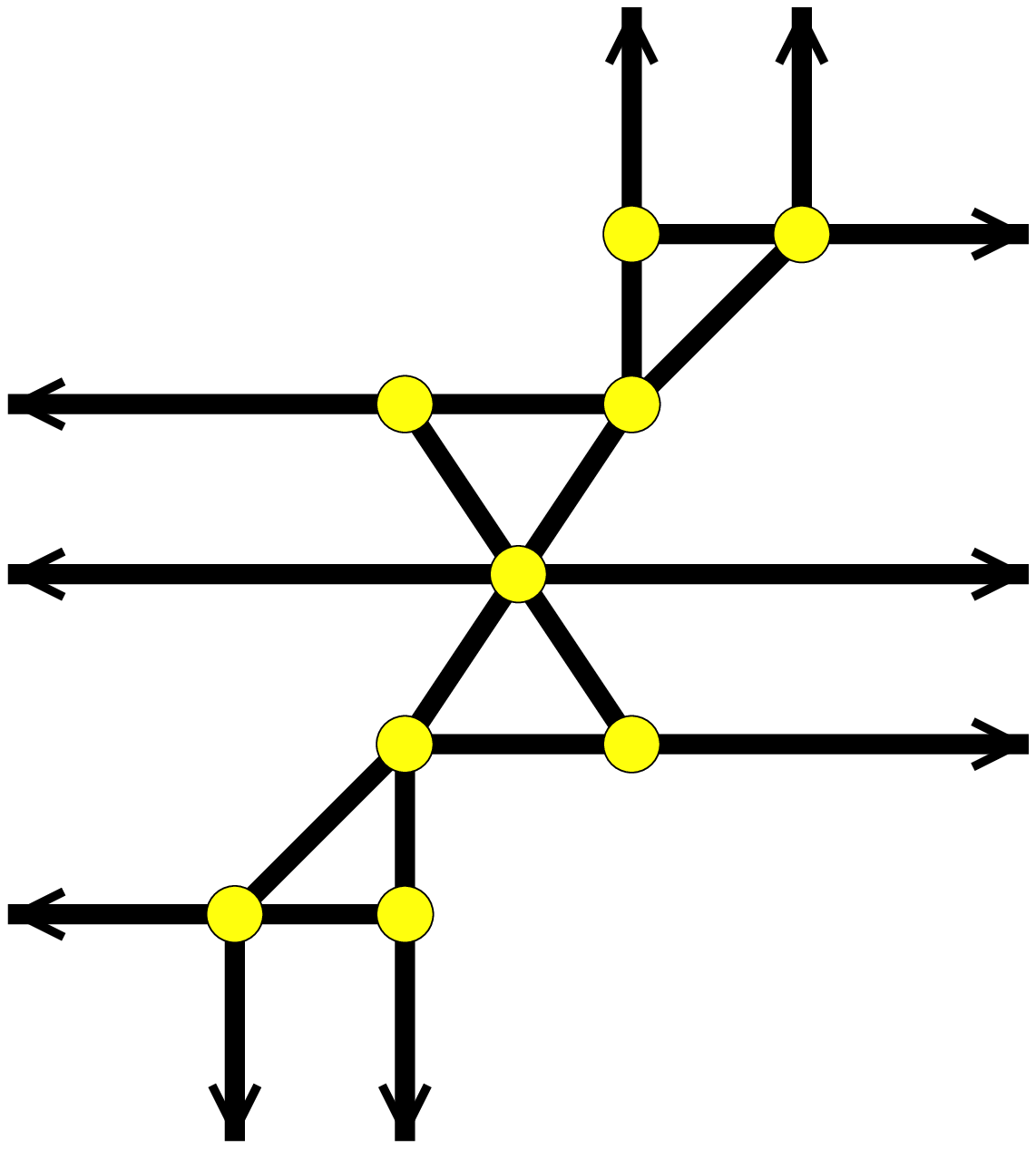,width=3cm}
    \hspace{1cm}
    \epsfig{figure=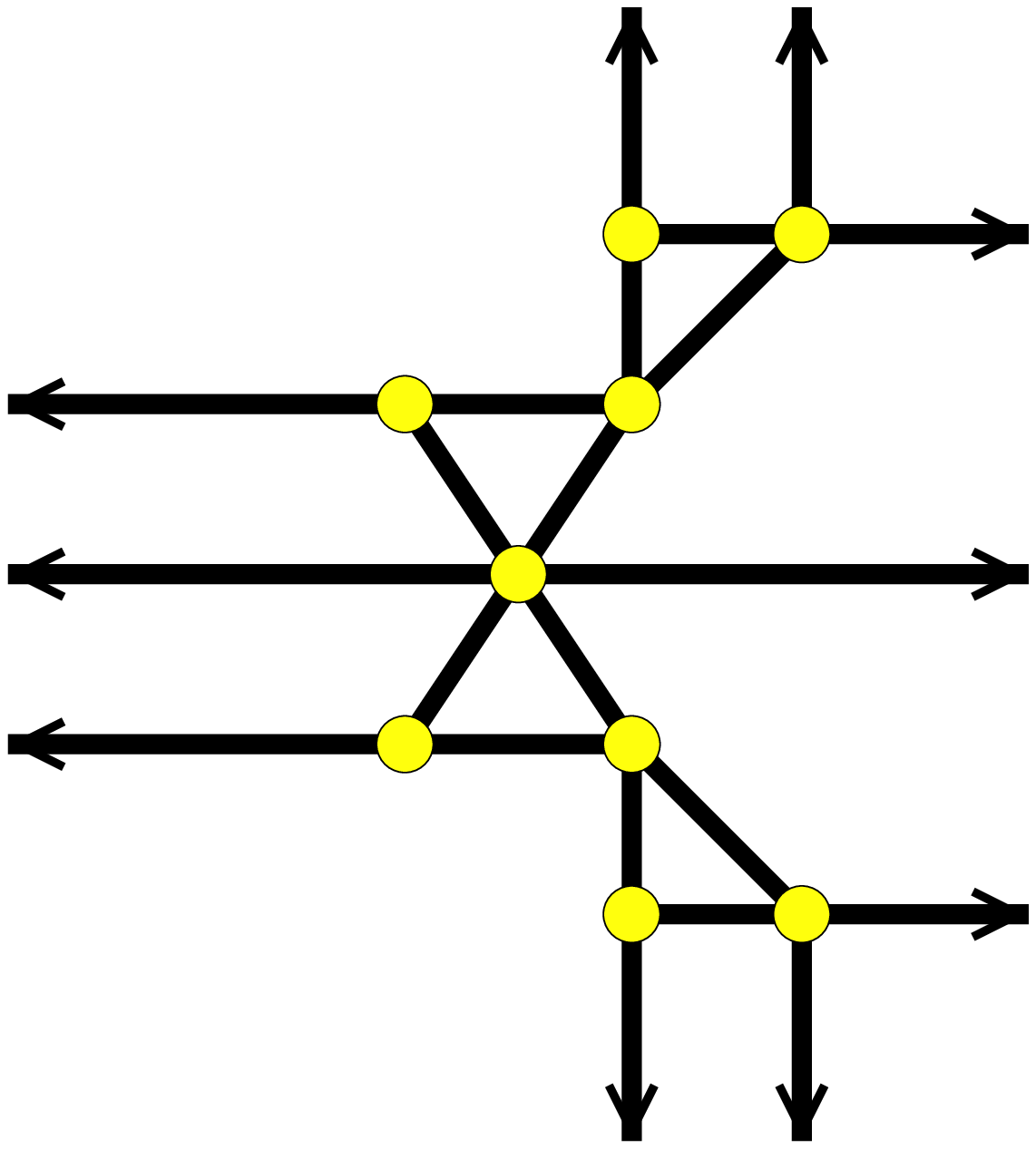,width=3cm}
    \caption{An example of two combinatorially different $3$-polyhedra
      with isomorphic vertex-facet incidences. The figures indicate the
      $1$-skeleta of the polyhedra.}
    \label{fig:3nonreconstr}
  \end{center}
\end{figure}

Propositions~\ref{prp:detectDimThree} and~\ref{prp:bndFacetDim} might
suggest to ask if the entire combinatorial structure of a
$d$-polyhedron can be reconstructed from its vertex-facet incidences
if $d=3$ or if $P$ has a bounded facet. However, the example given in
Figure~\ref{fig:3nonreconstr} shows that both answers are ``no''. 
The crucial feature of the example is that one can reflect the
``lower'' parts in the drawings without affecting the vertex-facet
incidences while changing the face poset (e.g., in contrast
to the left polyhedron the right one has two adjacent unbounded facets
that have three vertices each). For three-dimensional polyhedra this
is more or less the only kind of ambiguity that can arise.

\begin{prp}
  Given the vertex-facet incidences of a $3$-polyhedron $P$, for which
  $\fingr{P}$ is $2$-connected, one can determine $\facposet{P}$. 
\end{prp}

\begin{proof}
  One can compute $\fingr{P}$ from the vertex-facet incidences of $P$,
  and thus, one finds the graph of each facet of $P$. If all these
  graphs of facets are cycles then $P$ is bounded and the statement is
  clear. Otherwise, due to the $2$-connectedness of $\fingr{P}$, there
  is a unique (up to reorientation) way to arrange the paths that are
  the graphs of the unbounded facets of $P$ as a cycle. From this
  cycle, it is easy to determine the incidences of extremal rays and
  facets of $P$, which then allow to reconstruct the entire
  combinatorial structure of $P$.
\end{proof}

\begin{figure}[ht]
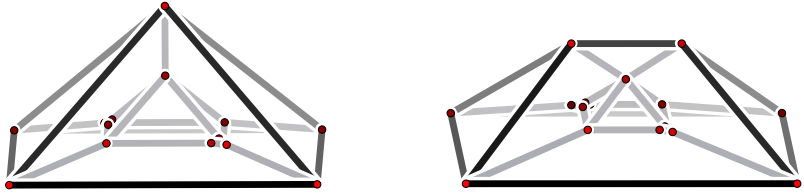

  \begin{center}
    \leavevmode
    \epsfig{figure=poly-1-bnd.eps,width=5cm}\hspace{1cm}
    \epsfig{figure=poly-2-bnd.eps,width=5cm}
    \caption{Schlegel diagrams (produced using
      \polymake~\cite{pm00,GJ00} and \javaview~\cite{PKPR00,Pol00}) illustrating two
      4-polyhedra $P_1$ and $P_2$ that have the same vertex-facet
      incidences, but different face posets.}
  \label{fig:sec-1:schlegels}
  \end{center}
\end{figure}

In larger dimensions, however, it is not true that higher
connectedness of the graph of a polyhedron is a sufficient condition
for the possibility to reconstruct its combinatorial structure from
its vertex-facet incidences. Figure~\ref{fig:sec-1:schlegels} shows
Schlegel-diagrams of (truncations of) two unbounded $4$-polyhedra.
These two polyhedra have the same vertex-facet incidences and a
$3$-connected graph, although their face posets are different (e.g.,
the right polyhedron has an extremal ray more than the left
one)\footnote{The data of these polyhedra as well as explanations on
  their construction can be found in the EG-Models archive at:
  \texttt{http://www-sfb288.math.tu-berlin.de/\allowbreak
    eg-models/\allowbreak models/\allowbreak polytopes/\allowbreak
    2000.05.001/\allowbreak \_preview.html}}.

The examples illustrated in Figures~\ref{fig:3nonreconstr}
and~\ref{fig:sec-1:schlegels} show that ``cone-like'' polyhedra are
not the only ones that cannot be reconstructed from their vertex-facet
incidences (even not in dimensions three and four). The polyhedra in
both examples are quite different from cones; each of them has a
full-dimensional vertex set, bounded facets, and the property that no
two facets have the same vertex set. Furthermore, in the
four-dimensional example, the vertex sets of the facets
even form an anti-chain (as promised in the introduction).

Nevertheless, any ambiguities in reconstructing the face poset of an
unbounded polyhedron from its vertex-facet incidences arise from some
degeneracy of~$P$.

\begin{thm}
  Given the vertex-facet incidences of a simple polyhedron $P$, one can 
  determine $\facposet{P}$.
\end{thm}

\begin{proof}
  Let $v$ be a vertex of a simple $d$-polyhedron $P$ and let
  $F_1,\dots,F_d$ be the facets of $P$ that contain $v$. Then
  the edges and extremal rays containing $v$ are precisely
  $$
  \bigcap_{i\in\{1,\dots,d\}\setminus\{i_0\}}F_i\qquad
         (i_0=1,\dots,d)\enspace.
  $$
  Since we can compute the edges of $P$ from a
  vertex-facet incidence matrix, we can thus also deduce
  (combinatorially) the extremal rays of $P$ and the information which
  ray is contained in which facets. From that, we can deduce
  the entire face poset of $P$. 
\end{proof}

Again, the example of cones shows that without dimension information
one can (in general) not decide from the vertex-facet incidences of a
polyhedron if it is simple.

All algorithms described in this section can be implemented 
such that their running time is bounded by a polynomial in
$|\facvertposet{P}|$.

To summarize the results in this section: we presented large classes
of (unbounded) polyhedra for which the combinatorial structures can be
reconstructed from their vertex-facet incidences as well as several
examples of polyhedra, for which this is not possible. Unfortunately,
these results do not yield a \emph{characterization} of the class of
those polyhedra that allow such reconstructions.


%% file: unbounded.tex
\section{Detecting Boundedness}\label{sec:unbounded}

In this section, we show that one can decide from the vertex-facet
incidences of a pointed polyhedron $P$ whether it is bounded or not.
It turns out that this only depends on the Euler characteristic of
(the order complex of) $\facvertposet{P}$. Thus, it can be read off
from the M\"obius function of $\facvertposet{P}$.

We recall some basic facts from topological combinatorics (see
Bj\"orner~\cite{Bj:TopMethods}).  Let $\Pi$ be a finite poset.  The
\emph{order complex}~$\Delta(\Pi)$ of~$\Pi$ is the finite simplicial
complex of all chains in~$\Pi$. We will use terminology from topology
in the context of finite posets such as~$\Pi$.  Throughout, this is
meant to refer to $\norm{\Delta(\Pi)}$ (i.e., any geometric
realization of $\Delta(\Pi)$, endowed with its standard topology).

It is well-known that the order complex~$\Delta(\facposet{P})$ of a
bounded $d$-polytope~$P$ is isomorphic (as a simplicial complex) to
the barycentric subdivision of the boundary~$\partial P$ of~$P$. In
particular, the topological type of $\facvertposet{P}$ is well-known
in this case.

\begin{lem}\label{lem:bounded}
  If $P$ is a $d$-polytope, then $\facposet{P}$ is homeomorphic to the
  $(d-1)$-sphere.
\end{lem}

If $P$ is an unbounded (pointed) polyhedron, then we can consider
$\facposet{P}$ as the sub-poset of $\facposet{\projclo{P}}$ that
consists of all faces of $\projclo{P}$ that are not contained in
$\farface$. Thus, we will identify $\Delta(\facposet{P})$ with the
sub-complex of $\Delta(\facposet{\projclo{P}})$ that is induced by all
chains $\{F\}$, where $F$ is a face of $\projclo{P}$ with
$F\not\subseteq\farface$.

\begin{lem}\label{lem:unbounded}
  If $P$ is an unbounded (pointed) polyhedron, then
  $\facposet{P}$ is contractible.
\end{lem}

\begin{proof}
  By Lemma~\ref{lem:bounded}, $\norm{\Delta(\facposetprojclo{P})}$ is
  homeomorphic to a sphere.  The induced subcomplexes
  $A=\Delta(\facposet{P})$ and $B=\Delta(\facposet{\farface})$ cover
  all vertices (i.e., one-element chains of $\facposet{\projclo{P}}$)
  of $\Delta(\facposet{\projclo{P}})$.  Using barycentric coordinates,
  it is seen that
  $\norm{\Delta(\facposet{\projclo{P}})}\setminus\norm{B}$ retracts
  onto $\norm{A}$. Thus, $\norm{A}$ has the same homotopy-type as
  $\norm{\Delta(\facposet{\projclo{P}})}\setminus\norm{B}$, where the
  latter is a simplicial sphere minus an induced ball.  Hence,
  $\facposet{P}$ is contractible.
\end{proof}

The two lemmas allow one to distinguish between the face
posets of bounded and unbounded polyhedra. Of course, there are
simpler ways to decide whether a face poset belongs to a bounded or
to an unbounded polyhedron (e.g., checking if every rank one element
is a join). However, in general we cannot reconstruct the face poset
of a polyhedron $P$ from its vertex-facet incidences (see
Section~\ref{sec:reconstruction}). Instead, we need criteria allowing
to distinguish between bounded and unbounded polyhedra that can be
computed from ~$\facvertposet{P}$. It turns out that the topological
criteria provided by Lemmas~\ref{lem:bounded} and~\ref{lem:unbounded}
can be exploited for this. 

Consider the poset maps $\phi:\facposet{P}\to\facvertposet{P}$,
mapping a face $F$ of a pointed polyhedron~$P$ to $\verts{F}$, and
$\psi:\facvertposet{P}\to\facposet{P}$, mapping the vertex set $S$ of
a face to the minimal face containing~$S$. Both $\phi$ and $\psi$ are
order preserving.  Moreover, $\phi(\psi(S))=S$ and $\psi(\phi(F))\leq
F$.

\begin{lem}\label{lem:homotopic}
  Let~$P$ be a pointed polyhedron. Then the face poset~$\facposet{P}$
  is homotopy-equivalent to the poset~$\facvertposet{P}$.
\end{lem}

\begin{proof}
  Setting $f(F)=\psi(\phi(F))$ defines an order preserving map from
  $\facposet{P}$ into itself such that each face~$F$ is comparable
  with its image~$f(F)$.  From the Order Homotopy
  Theorem~\cite[Corollary~10.12]{Bj:TopMethods}, we infer that
  $\facposet{P}$ is homotopy-equivalent to the
  image~$f(\facposet{P})$.  In fact, $f(f(F))=f(F)$, and hence
  $f(\facposet{P})$ is a strong deformation retract of~$\facposet{P}$.
  This proves the lemma, since $\psi$ is a poset isomorphism
  from~$\facvertposet{P}$ onto
  $\psi(\facvertposet{P})=f(\facposet{P})$.
\end{proof}

The reduced Euler characteristic of (the order complex of) a
poset~$\Pi$ is denoted by $\widetilde{\chi}(\Pi)$, i.e., 
$$
\widetilde{\chi}(\Pi)=\sum_{i=-1}^{D}(-1)^if_i(\Delta(\Pi))
$$
(where $f_i(\Delta(\Pi))$ is the number of $i$-faces
of~$\Delta(\Pi)$, and~$D$ is the dimension of~$\Delta(\Pi)$).  The
following result in particular shows that a polytope and an unbounded
polyhedron cannot have isomorphic vertex-facet incidences.

\begin{thm}\label{thm:decide}
  Let $P$ be a pointed polyhedron. Then $P$ is bounded if
  and only if $\widetilde{\chi}(\facvertposet{P})\not= 0$.
\end{thm}

\begin{proof}
  The reduced Euler characteristic of a $(d-1)$-sphere equals
  $(-1)^{d-1}$, while the reduced Euler characteristic of a
  contractible space vanishes.
  Thus the claim follows from Lemma~\ref{lem:bounded},
  Lemma~\ref{lem:unbounded}, and Lemma~\ref{lem:homotopic}.
\end{proof}

As an example consider the case where the unbounded polyhedron~$P$ has
a face~$F$ which contains all vertices of $P$.  Then
$\Delta(\facvertposet{P})$ is a cone over~$F$ (in the sense of
simplicial topology); in particular, it is
contractible and thus $\widetilde{\chi}(\facvertposet{P})=0$. 

The reduced Euler characteristic of the poset $\facvertposet{P}$ can
be computed efficiently as follows.  By adjoining an artificial top
element~$\hat{1}$ and an artificial bottom element~$\hat{0}$, the
poset~$\facvertposet{P}$ becomes a lattice~$\facvertlat{P}$.  Note
that we adjoin~$\hat{1}$ also in the case where $\facvertposet{P}$
already has a top element corresponding to a face
containing all  vertices of~$P$.

For every element $S\in\facvertlat{P}$ we define the~\emph{M\"obius
  function}, see Rota~\cite{Rot64} and Stanley~\cite{Stanley},
$$
\mu(S)=
\begin{cases}
  1 & \text{if }S=\hat{0}\enspace,\\
  -\displaystyle\sum_{S'\subsetneq S}\mu(S') 
  & \text{otherwise}\enspace.
\end{cases}
$$
The \emph{M\"obius number} $\mu(\facvertposet{P})=\mu(\hat{1})$ of
$\facvertposet{P}$ can be computed in time bounded polynomially in
$|\facvertposet{P}|$.  Since it is well-known (see
Stanley~\cite[3.8.6]{Stanley}) that
\begin{equation}
  \label{eq:muEuler}
  \mu(\facvertposet{P})=\widetilde{\chi}(\facvertposet{P})\enspace,
\end{equation}
this proves the following complexity result.

\begin{cor}
  There is an algorithm that decides for every vertex-facet incidence
  matrix of a polyhedron $P$ if $P$ is bounded. Its running time is
  bounded by a polynomial in $|\facvertposet{P}|$.
\end{cor}

Actually, Theorem~\ref{thm:decide} allows to decide even more from the
vertex-facet incidences of a polyhedron $P$.  Once we have computed
$\facvertposet{P}$ we clearly can also determine $\facvertlat{F}$ for
every facet $F$ of $P$ (since we know $\verts{F}$ for every facet $F$
of $P$). This is the interval between $\hat{0}$ and $\verts{F}$ in the
lattice $\facvertlat{P}$, where we have to add an additional top
element $\hat{1}$ if there is some other facet $F'$ of $P$ containing
$\verts{F}$.

\begin{cor}
  \label{cor:facetDecide}
  There is an algorithm that tells from a vertex-facet incidence
  matrix of a polyhedron~$P$ which facets of $P$ are bounded. Its
  running time is bounded by a polynomial in $|\facvertposet{P}|$.
\end{cor}


%% file: simpsimp.tex
\section{Simple and Simplicial Polyhedra}
\label{sec:simpsimp}

It is a well-known fact \cite[Exerc.~0.1]{GMZ} that a
$d$-\emph{polytope} which is both simple and simplicial is a simplex
or a polygon. Both properties (simplicity as well as simpliciality)
can be viewed as properties of vertex-facet incidences (see
Section~\ref{sec:definitions}). In this section, we generalize the
known result on polytopes to not necessarily bounded $d$-polyhedra
with $d\geq 2$.

\begin{thm}\label{thm:simple_simplicial}
  For $d \geq 2$, every simple and simplicial $d$-polyhedron is a
  simplex or a polygon.  In other words, unbounded simple and
  simplicial polyhedra do not exist.
\end{thm}

Our proof of Theorem \ref{thm:simple_simplicial} is organized into two
parts. The first part shows that the graph $\Gamma_P$ of a simple and
simplicial polyhedron $P$ is either a complete graph or a cycle.  In
the second part, we further deduce that a simple and simplicial
polyhedron has a \emph{circulant} vertex-facet incidence matrix. The
proof of Theorem~\ref{thm:simple_simplicial} is then completed by
showing that no unbounded $d$-polyhedron (with $d\geq 2$) can have a
circulant vertex-facet incidence matrix. Furthermore,
Propositions~\ref{prop:graph_cycle_or_complete} and~\ref{prop:circInc} 
yield characterizations of those polyhedra that have circulant
vertex-facet incidence matrices. 

\subsection{Graphs of Simple and Simplicial Polyhedra}

Throughout this section, let $P$ be a pointed simple and simplicial
$d$-polyhedron with $n$ vertices and $d\ge2$.  Double counting yields
that $P$ must also have $n$ facets. In particular, we have $n>d$
(since otherwise $P$ would be a cone, which is simple and simplicial
only for $d=1$). We denote by $V_P=\verts{P}$ the set of vertices of
$P$. For $S\subseteq V_P$ let $\inFacs{S}$ be the set of all
facets of $P$ that contain $S$. Recall that (since $P$ is simple) two
vertices $v$ and $w$ of $P$ form an edge if and only if
$\card{\inFacs{\{v,w\}}}=d-1$.

\begin{lem}\label{lem:no_equal_rows}
  Two different facets of $P$ cannot have the same set of vertices.
\end{lem}

\begin{proof}
  Suppose that there are two facets $F_1$ and $F_2$ of $P$ ($F_1\not=
  F_2$) with $\verts{F_1}=\verts{F_2}=:S$. Since $n > d$, and since
  $\Gamma_P$ is connected, there must be a vertex $v \notin S$ that is
  a neighbor of some vertex $w\in S$. Hence, we have
  $\card{\inFacs{\{v,w\}}}=d-1$. Because of $\card{\inFacs{\{w\}}}=d$
  and $F_1,F_2\in\inFacs{\{w\}}\supseteq\inFacs{\{v,w\}}$ this implies
  $F_1\in\inFacs{\{v,w\}}$ or $F_2\in\inFacs{\{v,w\}}$, which in both
  cases yields a contradiction to $v\not\in S$.
\end{proof}

For $S \subseteq V_P$, define $\touchedFacs{S}$ to be the set of those
facets of $P$ that have non-empty intersection with $S$. 

\begin{lem}\label{lem:span_of_ones}
  Let $S \subseteq V_P$ with $\card{S}>0$. Then
  $\card{\touchedFacs{S}} \geq \min \{ n,d+\card{S}-1\}$.
\end{lem}

\begin{proof}
  If $\card{\touchedFacs{S}} = n$, then the claim obviously is
  correct. Therefore, assume $\card{\touchedFacs{S}} < n$.  Since
  $\Gamma_P$ is connected, the vertices in $V_P \setminus S = \{z_1,
  \dots, z_{r}\}$ ($r= n-\card{S}$) can be ordered such that $z_{i+1}$
  is adjacent to some vertex of $S_i = S \cup \{ z_1, \dots, z_i \}$
  for each $i \in \{0, \dots, r-1\}$ (additionally, define $S_r = S
  \cup \{z_1, \dots, z_r\}$). Clearly $\card{\touchedFacs{S_i}} \leq
  \card{\touchedFacs{S_{i-1}}}+1$, since vertex $z_i$ has $d-1$ facets
  in common with some vertex in $S_{i-1}$.
  
  Define $l$ to be the last $i$, such that $\card{\touchedFacs{S_i}} =
  \card{\touchedFacs{S_{i-1}}}+1$, i.e., $l$ is the last index, where we
  encounter a new facet ($l$ is well-defined due to
  $\card{\touchedFacs{S}} < n$). Since this facet must contain $d-1$
  vertices from $V_P\setminus S_l$, we have $r-l \geq d-1$, which
  yields $n-l \geq d+\card{S}-1$.
  
  Furthermore, we have $\card{\touchedFacs{S}} + l \geq n$, since $S_l$
  intersects all facets. It follows $\card{\touchedFacs{S}} \geq
  n-l \geq d+\card{S}-1$.
\end{proof}

For $S\subseteq V_P$ let $\Gamma_P(S)$ be the subgraph of $\Gamma_P$
induced by $S$.

\begin{lem}\label{lem:full_column_sum}
  Let $S \subset V_P$ with $0<\card{S}\leq d$, such that
  $\Gamma_P(S)$ is connected. Then $\card{\inFacs{S}}=d-\card{S}+1$ holds.
\end{lem}

\begin{proof}
  Since $\Gamma_P(S)$ is a connected subgraph of the connected graph
  $\Gamma_P$ (which has $n>d$ vertices), there is a chain $\emptyset
  \subsetneq S_1 \subsetneq S_2 \subsetneq \ldots \subsetneq S_d$ with
  $S_{\card{S}}=S$, such that $\card{S_i}=i$ and $\Gamma_P(S_i)$ is connected
  for all~$i$.
  
  For every $1<i\leq d$, the vertex $v$ with $S_i\setminus
  S_{i-1}=\{v\}$ is connected to some vertex $w\in S_{i-1}$. From
  $\card{\inFacs{\{w\}}\setminus\inFacs{\{v\}}}=1$ we infer
  $\card{\inFacs{S_{i-1}}\setminus\inFacs{\{v\}}}\leq 1$, and thus,
  $\card{\inFacs{S_i}}\geq\card{\inFacs{S_{i-1}}}-1$.  Together with
  $\card{\inFacs{S_1}}=d$ (since $P$ is simple) and
  $\card{\inFacs{S_d}}\leq 1$ (by Lemma~\ref{lem:no_equal_rows}), this
  implies $\card{\inFacs{S_i}}=d-i+1$ for all $1\leq i\leq d$.
\end{proof}

The next three lemmas show that $\Gamma_P$ has a very special
structure.

\begin{lem}
  If $\fingr{P}$ contains a cycle $C$ of size $k>d$, then 
  $\fingr{P}$ is the cycle $C$ or a complete graph on $n=d+1$ nodes.
\end{lem}

\begin{proof}
  Let $C = ( v_0, \dots, v_{k-1}, v_0 )$ be a cycle of size $k > d$ in
  $\Gamma_P$.  In the following, all indices are taken
  modulo~k.  For $0 \le i \le k-1$ define the set $C_i = \{ v_i,
  \dots, v_{i+d-1}\}$ of size $d$.  Clearly, $\Gamma_P(C_i)$ is
  connected, and, by Lemma~\ref{lem:full_column_sum}, there exists
  exactly one facet $F_i$ with $\inFacs{C_i}=\{F_i\}$.  Due to $k >
  d$, the facets $F_0, \dots, F_{k-1}$ are pairwise distinct. This
  means that $\inFacs{\{v_i\}}=\{F_{i-d+1},\dots,F_i\}$ (since $P$ is
  simple) and $\verts{F_i}=C_i$ (since $P$ is simplicial).  Hence,
  every vertex that is adjacent to one of the nodes $v_0,\dots,
  v_{k-1}$ must be contained in at least one (more precisely, in
  $d-1>0$) of the facets $F_0,\dots, F_{k-1}$, and thus it lies in
  $\{v_0, \dots, v_{k-1}\}$.
  
  Since $\Gamma_P$ is connected, this means that $n = k$.  For $n=d+1$
  this immediately yields that $\Gamma_P$ is a complete graph on
  $n=d+1$ nodes, while for $n>d+1$ one finds that $\Gamma_P$ is the
  cycle~$C$ (since, in this case,
  $\card{\inFacs{\{v_i\}}\cap\inFacs{\{v_j\}}}=d-1$ if and only if
  $j\equiv i\pm 1\mod k$).
\end{proof}

\begin{lem}
  If $\fingr{P}$ contains a cycle of length $k\leq d$, then
  $\fingr{P}$ is a complete graph on $n=d+1$ nodes.
\end{lem}

\begin{proof}
  Let $\widetilde{C} = (v_0, \dots, v_{k-1}, v_0)$ be a cycle in $\fingr{P}$ of
  size $k\leq d$. For each $i\in\{0,\dots,k-1\}$ define
  $\widetilde{C}_i=\{v_0,\dots,v_i\}$. Taking all indices modulo~$k$, we
  have $|\inFacs{\{v_i,v_{i+1}\}}|=d-1$ for each $i$, and hence, there
  are facets $F_i$ and $G_i$ with 
  $$
  \inFacs{\{v_i\}}\setminus\inFacs{\{v_{i+1}\}}=\{F_i\}
  \quad\text{and}\quad
  \inFacs{\{v_{i+1}\}}\setminus\inFacs{\{v_i\}}=\{G_i\}\enspace.
  $$
  It follows that
  \begin{equation}
    \label{eq:Ck-1}
    \inFacs{\widetilde{C}_{k-1}}=\inFacs{\widetilde{C}_0}\setminus\{F_0,\dots,F_{k-1}\}\enspace.    
  \end{equation}
  If $\fingr{P}$ is not complete, then $n>d+1$ holds, and we infer
  from Lemma~\ref{lem:span_of_ones} that $\card{\touchedFacs{\widetilde{C}_2}}\geq
  d+2$, which implies $G_0,G_1\not\in\inFacs{\widetilde{C}_0}$ (with $G_0\not=
  G_1$). Due to $\{F_0,\dots,F_{k-1}\}=\{G_0,\dots,G_{k-1}\}$,
  Equation~\eqref{eq:Ck-1} implies  
  $$
  \card{\inFacs{\widetilde{C}_{k-1}}}\geq\card{\inFacs{\widetilde{C}_0}}-(k-2)=d-k+2\enspace,
  $$
  contradicting Lemma~\ref{lem:full_column_sum}.
\end{proof}

By the above two lemmas, $\fingr{P}$ cannot contain any cycles, unless 
it is complete or a cycle itself. Thus, we are left with the case of
$\fingr{P}$ not containing any cycles at all. 

\begin{lem}\label{lem:graph_not_tree}
  $\Gamma_P$ is not a tree.
\end{lem}

\begin{proof}
  Assume $\Gamma_P$ is a tree. Let $v \in V_P$ be a leaf of
  $\Gamma_P$ with $u$ being the unique vertex of $\fingr{P}$ adjacent
  to $v$. Due to
  $\card{\inFacs{\{v\}}\setminus\inFacs{\{v,u\}}}=1$,
  there is one facet that induces a subgraph of $\fingr{P}$ in
  which $v$ is isolated.  This, however, is a contradiction to
  Lemma~\ref{lem:gamma_connected}.
\end{proof}

Altogether this proves the following.

\begin{prp}\label{prop:graph_cycle_or_complete}
  Let $P$ be a simple and simplicial $d$-polyhedron ($d \geq 2$) with
  $n$ vertices.  Then $\Gamma_P$  is an $n$-cycle or a complete
  graph on $n=d+1$ nodes.
\end{prp}

It is worth to mention that one can generalize Proposition
\ref{prop:graph_cycle_or_complete} in the following way. Let $A$ be a
$0/1$-matrix of size $n \times n$ with row and column sums~$d$. Define
a graph $\Gamma_A$ on the columns of $A$, such that two columns are
adjacent if and only if they have exactly $d-1$ ones in common rows.
Then, by the same arguments as above, one can show that the
connectedness of $\Gamma_A$ already implies that $\Gamma_A$ is a cycle
or a complete graph. The only difference in the proof arises in Lemma
\ref{lem:graph_not_tree}. Here one has to prove additionally that, for
each row, the subgraph of $\Gamma_A$ that is induced by the ones in
that row is connected (if $\Gamma_A$ is connected).

\subsection{Circulant Matrices}

We will now exploit Proposition \ref{prop:graph_cycle_or_complete} 
to show that every simple and simplicial polyhedron has a very special 
vertex-facet incidence matrix. 

Let $n,d$ be integers satisfying $1 \leq d \leq n$. The
$(n,d)$-\emph{circulant}~$M(n,d)$ is the $n\times n$-matrix with $0/1$
entries whose coefficients $m_{ij}$ ($i,j\in\{0,\dots,n-1\}$) are
defined as follows:

$$
m_{ij}= 
\begin{cases}
  1 & \text{if $j\in\{i,i+1\bmod n,\ldots,i+d-1\bmod n\}$}\\
  0 & \text{otherwise}
\end{cases}
$$

For $d \geq 1$, the $(d+1,d)$-circulant is an incidence matrix of the
$d$-simplex, and for $n \ge 3$, the $(n,2)$-circulant is an incidence
matrix of the ($2$-dimensional) $n$-gon.

\begin{prp}\label{prop:simpsimpcirc}
  A polyhedron~$P$ is simple and simplicial if and only if it has a
  circulant $M(n,d)$ as  a vertex-facet incidence matrix.  In this case,
  $\dim(P)=d$.
\end{prp}

\begin{proof}
For the ``if''-direction of the proof, let $P$ be a polyhedron with a
vertex-facet incidence matrix $M(n,d)$ ($1\leq d\leq n$). The cases
$d=1$ (implying $n\in\{1,2\}$) as well as $d=n$ (implying $d=n=1$) are
trivial. Therefore, let $2\leq d <n$. Obviously, it suffices to show
$\dim(P)=d$. To each row $i \in \{0,\dots,n-1\}$ of $M(n,d)$ there
corresponds a facet $F_i$ of $P$. For $0 \leq j \leq d-1$ define $G_j =
F_0 \cap \dots \cap F_j$. Clearly, $G_j \supseteq G_{j+1}$ holds for $0
\leq j < d-1$.  Due to $\verts{G_{j}} = \verts{F_0} \cap \dots \cap
\verts{F_j}$ it follows $\verts{G_j} \supsetneq \verts{G_{j+1}}$ and
therefore $G_{j} \supsetneq G_{j+1}$. Now $F_0 = G_0 \supsetneq G_1
\supsetneq \dots \supsetneq G_{d-1}$ is a (decreasing) chain of length
$d-1$ in the face poset of $P$. Hence we have $\dim P \geq d$.  Since
each vertex must be contained in at least $\dim P$ facets it follows
that $\dim P \leq d$ (because each vertex of~$P$ is contained in
precisely $d$~facets).
  
Conversely, let $P$ be a simple and simplicial $d$-polyhedron ($d\geq
1$) with $n$ vertices. The case $d=1$ is checked easily. Thus, assume
$d\geq 2$. By Proposition~\ref{prop:graph_cycle_or_complete},
$\fingr{P}$ either is a complete graph on $n=d+1$ nodes or it is a
cycle. In the first case, every vertex-facet incidence matrix of $P$
is the complement of a permutation matrix, which can be transformed to
$M(n,d)$ by a suitable permutation of its rows.  In the second case,
consider any vertex-facet incidence matrix $A$ of $P$, where the
columns are assumed to be ordered according to the cycle $\Gamma_P$.
Call two positions $(i,j)$ and $(i,k)$ in $A=(a_{fv})$
($f,v\in\{0,\dots,n-1\}$) \emph{mates} if $k \equiv j+1\pmod n$ and
$a_{ij}=a_{ik}=1$. Walking around the cycle $\Gamma_P$, we find that
the total number of mates in $A$ is precisely $n(d-1)$ (because every
edge is contained in precisely $d-1$ facets).  But then, since every
row of $A$ has only $d$ ones (because $P$ is simplicial), it follows
that in each row the ones must appear consecutively (modulo $n$).
Denote by $s(i)$ the starting position of the block of ones in
row~$i$. Because there are no equal rows in $A$ (by Lemma
\ref{lem:no_equal_rows}) we deduce that $s$ defines a permutation of
the rows of $A$ which tells us how to transform $A$ to $M(n,d)$.
\end{proof}

The following result finishes the proof of
Theorem~\ref{thm:simple_simplicial} 
(via Proposition~\ref{prop:simpsimpcirc}).

\begin{prp}
\label{prop:circInc}
  If a polyhedron~$P$ has $M(n,d)$ ($2\leq d < n$) as a vertex-facet
  incidence matrix, then $n=d+1$  ($P$ is a $d$-simplex) or $d=2$ 
  ($P$ is an $n$-gon).
\end{prp}

\begin{proof}
  If $n=d+1$, then $M(n,d)$ is a vertex-facet incidence matrix
  of a $d$-simplex. Hence, by Theorem~\ref{thm:decide}, $P$ cannot be
  unbounded, and thus it must be a $d$-simplex as well. Therefore, in the
  following we will assume $n>d+1$. 
  
  Let us first treat the case $d+1<n<2d-1$.  Consider the facets $F$
  and~$F'$ corresponding to rows~$0$ and~$n-d+1$, respectively.  If we
  identify the vertices of $P$ with the column indices
  $\{0,\dots,n-1\}$ of $M(n,d)$, then the vertex set of the face
  $G=F\cap F'$ is $\{0\}\cup\{n-d+1,\ldots,d-1\}$, where $\{n-d+1,
  \dots, d-1\} \neq \varnothing$ (due to $n<2d-1$). By
  Propositions~\ref{prop:simpsimpcirc}
  and~\ref{prop:graph_cycle_or_complete}, $\fingr{P}$ is an $n$-cycle
  (due to $n>d+1$). Since neither vertex~$1$ nor vertex~$n-1$, which
  are the only neighbors of~$0$ in $\fingr{P}$, are contained in~$G$,
  we conclude that the subgraph of $\fingr{P}$ induced by~$G$ is
  disconnected, which is a contradiction to
  Lemma~\ref{lem:gamma_connected}.

  Hence, we can assume $n\geq 2d-1$. This implies
  $$
  \facvertposet{P}=\SetOf{\{i,\dots,i+s-1\}}{i\in\{0,\dots,n-1\}, s\in\{1,\dots,d\}}
  $$
  (where, again, all indices are to be taken modulo~n), i.e.,
  $\facvertposet{P}$ consists of all (cyclic) intervals of
  $\{0,\dots,n-1\}$ with at least one and at most $d$ elements. We
  will compute the M\"obius function $\mu$ (see
  Section~\ref{sec:unbounded}) on the lattice $\facvertlat{P}$
  (which arises by adding artificial top and bottom elements $\hat{1}$
  and $\hat{0}$ to $\facvertposet{P}$). For each
  $s\in\{1,\dots,d\}$ let $\mu(s)=\mu(\{0,\dots,s-1\})$. Obviously,
  for every $F\in\facvertposet{P}$ with $\card{F}=s$ we have
  $\mu(F)=\mu(s)$.  In particular, one readily deduces $\mu(1)=-1$ and
  $\mu(2)=-(1+2\cdot(-1))=1$. For $3\leq s\leq d$ we then infer (by
  induction) $\mu(s)=-(1+s\cdot(-1)+(s-1)\cdot(+1))=0$. Thus, we finally
  calculate
  $$
  \mu(\facvertposet{P})=\mu(\hat{1})=-(1+n\cdot(-1)+n\cdot(+1))=-1\enspace,
  $$
  which by~\eqref{eq:muEuler} and Theorem~\ref{thm:decide} implies
  that $P$ is bounded (and, hence, an $n$-gon). 

  (Alternatively, one could derive from the \emph{Nerve Lemma}
  \cite[Theorem~10.7]{Bj:TopMethods} that $\facvertposet{P}$ is homotopy-equivalent
  to a circle for $n\geq 2d-1$, and thus, $P$ must be a polygon.)
\end{proof}
